\documentclass[a4paper,10pt]{article}
\usepackage[latin1]{inputenc}
\usepackage[T1]{fontenc}
\usepackage[french,english]{babel}

\usepackage{amsfonts}
\usepackage{amssymb}
\usepackage{eufrak}
\usepackage{amsmath}
\usepackage{enumerate}
\usepackage{theorem}
\usepackage[dvips]{graphicx}

\newtheorem{Def}{Définition}[section]
\newtheorem{Prop}[Def]{Proposition}
\newtheorem{Th}[Def]{Théorème}

\newtheorem{Lem}[Def]{Lemme}

\theoremstyle{plain}
\theorembodyfont{\itshape}

\theorembodyfont{\small\upshape}

\theoremstyle{plain}

\newcommand{\xe}[1]{\mathbf{#1}}

\DeclareMathOperator{\mult}{mult}

\title{Constantes de Seshadri du diviseur anticanonique des surfaces de del Pezzo}
\author{Amaël BROUSTET}

\begin{document}
\selectlanguage{french}
\maketitle
{\let\thefootnote\relax
\footnote{\hskip-0.74em
\textbf{Key-words :} Seshadri constants, rational curves.

\textbf{A.M.S.~classification :} 14J26, 14J45. 
 }}

\vspace{-1cm}
 
\begin{center}
\begin{minipage}{100mm}
\scriptsize
 
 {\bf Abstract.} 
Seshadri constants, introduced by Demailly, capture positivity of a nef divisor at a point. We compute in this note the Seshadri constants of the anticanonical bundle at every point of Del Pezzo surfaces. During the proof, we enlight the role of rational curves in our computations. We present then two exemples where the positivity of the anticanonical bundle cannot be detected using rational curves. 
 \end{minipage}
 \end{center}



\section{Introduction}

Le concept de positivité locale, introduit par J.P Demailly \cite{Dem1}, consiste à mesurer au travers des constantes de Seshadri la positivité d'un diviseur nef en un point. 

\begin{Def}
\label{Prop : inf sur les courbes}
Soit $X$ une variété projective lisse, $x$ un point de $X$ et $D$ un diviseur nef sur $X$. La constante de Seshadri en $x$ de $D$ est le réel positif :
$$\varepsilon (D;x) = \inf_{x\in C \subset X} \frac{D \cdot C}{\mult_x C}$$
la borne inférieure étant prise sur l'ensemble des courbes irréductibles $C \subset X$ passant par $x$.
\end{Def}

Pour plus d'informations sur les constantes de Seshadri, on recommande le chapitre $5$ de $\cite{Laz}$.

Si $x$ est un point en position générale, la constante ne dépend pas du point. Dans cette note, on calcule \emph{en tout point} les constantes de Seshadri du diviseur anticanonique des surfaces de del Pezzo lisses. On peut dans notre cas donner une définition précise de la notion de point en position générale.

\begin{Def}
Dans la suite, on notera $\mu_r : X_r \longrightarrow \mathbb{P}^2 $ l'éclatement du plan en $r$ points distincts $x_1, \ldots ,x_r$.  Un ensemble de $r\leqslant 8$ points $\{x_1,\ldots ,x_r\}$ du plan est dit en position générale si aucun sous-ensemble de $3$ de ces points n'est sur une droite, si aucun sous-ensemble de $6$ d'entre eux n'est sur une conique et si $8$ d'entre eux ne sont pas sur une cubique singulière en l'un deux.  Lorsque $r\leqslant 7$, un point $x$ de $X_r$ sera dit en position générale si son image $y$ par $\mu_r$ est distincte des points d'éclatement $x_i$ et si les points $\{x_1,\ldots ,x_r,y\}$ sont en position générale. Si le point $x$ n'est pas en position générale, il est alors sur une courbe exceptionnelle ou sur la transformée stricte d'une des courbes précédemment citée. On appelle cette courbe la \emph{courbe distinguée} contenant $x$.
\end{Def}

On montre le résultat suivant :

\begin{Th}
\label{Prop : résultat}
Si $r\leqslant 5$,
la constante de Seshadri de $-K_{X_r}$ au point $x$ vaut :
\begin{itemize}
\item[$\centerdot$] $\varepsilon(-K_{X_r};x)=2$  si $x$ est en position générale,
\item[$\centerdot$]  $\varepsilon(-K_{X_r};x)=1$ sinon.
\end{itemize}
Si $r=6$,
la constante de Seshadri de $-K_{X_6}$ au point $x$ vaut :
\nopagebreak
\begin{itemize}
\item[$\centerdot$] $\varepsilon(-K_{X_6};x)=3/2$ si $x$ est en position générale,
\item[$\centerdot$]  $\varepsilon(-K_{X_6};x)=1$ sinon.
\end{itemize}
Si $r=7$,
la constante de Seshadri de $-K_{X_7}$ au point $x$ vaut :
\nopagebreak
\begin{itemize}
\item[$\centerdot$] $\varepsilon(-K_{X_7};x)=4/3$ si $x$ est en position générale,
\item[$\centerdot$]  $\varepsilon(-K_{X_7};x)=1$ sinon.
\end{itemize}
Si $r=8$,
les constantes de Seshadri de $-K_{X_8}$ valent $\frac{1}{2}$ en au plus $12$ points n'appartenant pas au diviseur exceptionnel et $1$ partout ailleurs.
\end{Th}

Au cours de la preuve du théorème \ref{Prop : résultat}, nous mettrons en relief le rôle des courbes rationnelles dans le calcul des constantes de Seshadri. Plus exactement nous obtenons le résultat suivant :
\begin{Prop}
\label{Prop : CR}
Si $S$ est une surface de del Pezzo lisse
$$\varepsilon (-K_S;x)=\inf_{x\in C\textrm{ rationnelle}} \frac{-K_S \cdot C}{\mult_x C}.$$
\end{Prop}

Malheureusement, les courbes rationnelles ne permettent en général pas de détecter la positivité sur une variété rationnellement connexe. On donne à la fin de cette note deux exemples de surfaces rationnellement connexes dont le diviseur anticanonique possède une intersection positive avec toute courbe rationnelle mais n'est pas pas nef, ni même pseudoeffectif dans le cas du deuxième exemple. 

\section{Constantes de Seshadri du diviseur anticanonique des surfaces de Fano}
Le théorème suivant (\cite{Ful}, théorème 1 page 110) servira dans la suite :
\begin{Th}
\label{Th:Courbes_points_pos_generale}
Soit $x_1, \ldots, x_n$ $n$ points distincts de $\xe{P}^2$.
On note $V(d;r_1 x_1, \ldots , r_n x_n)$ l'espace projectif des courbes de degré $d$ et de multiplicité au moins $r_i$ en $x_i$.
C'est un sous-espace projectif de $\xe{P}^{\frac{d(d+3)}{2}}$ et
$$\dim V(d;r_1 x_1, \ldots , r_n x_n) \geqslant \frac{d(d+3)}{2} - \sum_{i=1}^n \frac{r_i(r_i +1)}{2}.$$
En particulier $V(d;r_1 x_1, \ldots , r_n x_n)$ n'est pas vide si $\frac{d(d+3)}{2} \geqslant \sum_{i=1}^n \frac{r_i(r_i +1)}{2}$.
\end{Th}

De même, on utilisera aussi la proposition suivante (\cite{Lnm},  théorème 1 p39).
\begin{Prop}
\label{8lis}
Soient $x, x_1, \ldots ,x_7$ huit points du plan tels que :
\begin{itemize}
\item quatre de ces points ne soient pas alignés,
\item aucune conique ne passe par $7$ d'entre eux.
\end{itemize}
Alors il existe une cubique lisse passant par ces huit points.
\end{Prop}

\subsection{Preuve du théorème \ref{Prop : résultat}}
\paragraph{Le cas $r\leqslant 6$.}

Le diviseur anticanonique $-K_{X_r}$ est très ample donc $\varepsilon (-K_{X_r};x)\geqslant 1$ pour tout point $x$.
Si $x$ n'est pas en position générale, la courbe distinguée $C$ contenant $x$ vérifie $-K_X \cdot C = 1$. On en déduit que $\varepsilon (-K_{X_r};x)= 1$

Si $x$ est en position générale, il existe alors un membre irréductible et réduit $D_x \in |-K_{X_r}|$ vérifiant $\mult_x D_x = 2$. Supposons ceci vrai pour le moment, on a alors pour toute courbe $C$ contenant $x$ différente du support de $D_x$ $$D_x \cdot C \geqslant 2 \mult_x C.$$
De plus $D_x^2=9-r\geqslant 4$ si $r\leqslant 5$ et  $D_x^2=6$ si $r=6$.
D'où  $\varepsilon (-K_{X_r};x)= 2$ si $r\leqslant 5$ et $\varepsilon (-K_{X_r};x)= \frac{3}{2}$ si $r=6$.

Il reste à montrer l'existence de $D_x$. D'après \ref{Th:Courbes_points_pos_generale}, il existe une courbe plane $D$ passant par tous les $x_i$ et dont la multiplicité au point $x$ est supérieure à $2$. Il suffit de vérifier que $\mult_{\mu_r(x)} D \leqslant 2$ et que $D$ est bien irréductible et réduite, $D_x$ sera alors la transformée stricte de $D$. Quitte à compléter l'ensemble des points $x_i$, on peut supposer $r=6$. Si $D$ n'est pas irréductible et réduite, $D$ est l'union de trois droites ou d'une droite et d'une conique. Par la position générique de l'ensemble de points $\{x_1,\ldots,x_6,x\}$, la courbe $D$ ne peut alors passer par tous les points avec la multiplicité prescrite.

La courbe $D$ étant irréductible, son intersection avec une droite $L$ passant par $\mu_r(x)$ et un autre point $z\in D$ vaut au moins
$$D\cdot L =\mult_z D + \mult_{\mu_r(x)} D$$
et au plus $3$. On en déduit que $\mult_{\mu_r(x)} D \leqslant 2$ et que $z$ est un point lisse de $D$.

\paragraph{Le cas $r=7$}
Pour tout point $x\in X_7$ , la proposition \ref{8lis} implique l'existence d'un membre $D_x \in | -K_{X_7}|$ passant par $x$ et lisse au point $x$, d'où $\varepsilon (-K_{X_r};x)\geqslant 1$. Si le point $x$ n'est pas en position générale, on déduit comme précédemment que $\varepsilon (-K_{X_r};x)= 1$.

Si le point $x$ est en position générale, il existe un membre irréductible et réduit $D_x\in |-2K_x|$ dont la multiplicité en $x$ vaut $3$. Cela revient à prouver l'existence d'une courbe plane de degré $6$ irréductible, réduite de multiplicité $3$ en $\mu_7(x)$ et de multiplicité $2$ au points $x_i$. Le théorème \ref{Th:Courbes_points_pos_generale} implique l'existence d'une courbe plane $D$ de degré $6$ et de multiplicités au moins égales à celles attendues aux points $x_i$ et $\mu_7(x)$. Si cette courbe est irréductible et réduite, son intersection avec une cubique passant par les $8$ points $x_1,\ldots ,x_7,\mu_7(x)$ et un neuvième sur la courbe $D$ vaut $18$ d'après le théorème de Bézout. On en déduit que $\mult_{x_i} D = 2$ pour tout $i$ et que $\mult_{\mu_r(x)} D = 3$. Il reste à montrer que cette courbe est irréductible et réduite. Si la courbe $D$ était l'union de deux cubiques $C_1$ et $C_2$ (éventuellement non irréductibles ou non réduites), les points $x_1,\ldots ,x_7$ étant en position générale, on aurait alors
$$17 = \sum_{i\leqslant 7} \mult_{x_i} D + \mult_{\mu_7(x)} D = \sum_{j=1,2} \,\sum_{i\leqslant 7} \mult_{x_i} C_j + \mult_{\mu_7(x)} C_j \leqslant 2\times 8.$$
Si $D$ n'est pas irréductible, $D$ est donc soit l'union de 3 coniques soit l'union d'une quintique réduite irréductible et d'une droite ou d'une quartique réduite irréductible et d'une conique. En intersectant cette quintique $Q$ avec une cubique $C$ passant par tous les points $x_1,\ldots ,x_7,x$ et un autre point de la quintique, on obtient d'après le théorème de Bézout
$$\sum_{i\leqslant 7} \mult_{x_i} Q + \mult_{\mu_7(x)} Q \leqslant 15 - 1 = 14.$$
Encore une fois, on aurait dans ce cas $$\sum_{i\leqslant 7} \mult_{x_i} D + \mult_{\mu_7(x)} D \leqslant 16$$
ce qui n'est pas possible. On procède de même pour les cas où $D$ est l'union d'une quartique et d'une conique et où $D$ est l'union de $3$ coniques. Le diviseur $D_x$ est donc irréductible et réduit. Comme au paragraphe précédent, on en conclut que  $\varepsilon (-K_{X_r};x)= \min \{\frac{3}{2},\frac{D_x^2}{3}\}=\frac{4}{3}$. De même qu'au paragraphe précédent, la courbe $D$ qui permet ``d'atteindre'' la constante de Seshadri est rationnelle.

\paragraph{Le cas $r=8$.}

Le système linéaire $|-K_{X_8}|$ est sans point base, de plus ses membres sont irréductibles et réduits grâce à la position générale des points $x_1,\ldots,x_6$. On en déduit que  $\varepsilon (-K_{X_r};x)= 1$ sauf aux éventuels points singuliers des membres de $|-K_{X_8}|$. Mais d'après la position générale des points $x_1,\ldots,x_6$ ces points singuliers sont en dehors du diviseur exceptionnel de $\mu_r$ et correspondent aux singularités des cubiques du pinceau de cubiques passant par les $x_i$. Le nombre de cubiques singulières dans un pinceau général de cubique est un problème classique de géométrie énumérative et vaut $12$. En ces points la constante de Seshadri de $-K_{X_8}$ vaut $\frac{1}{2}$.

Contrairement aux cas précédents, il n'existe en général pas pour les constantes de Seshadri de $-K_{X_8}$ de courbe rationnelle $\Gamma$ telle que $$ \varepsilon (-K_{X_8};x) = \frac{-K_{X_8} \cdot \Gamma}{\mult_x \Gamma}.$$
Pour voir cela, on peut notamment utiliser le lemme \ref{GP}.

On peut cependant noter qu'il existe une suite $(\Gamma_k)$ de courbes rationnelles telles que $$ \varepsilon (-K_{X_8};x) = \lim_{k\rightarrow \infty} \frac{-K_{X_8} \cdot \Gamma_k}{\mult_x \Gamma_k}.$$
C'est une conséquence directe de \cite{GLS}, lemme 3.2.10 :
\begin{Lem}
Soit $\bar{X}$ l'éclatement de $\xe{P^2}$ en $9$ points en position générale, on note $E_0$ le tiré en arrière de $\mathcal{O}_{\xe{P}_2}(1)$, et $E_i$ la préimage du point d'éclatement $x_i$. 
Pour tout entier $m \geq 1$ il existe une courbe rationnelle nodale irréductible dans le système
$$|3mE_0 - mE_1 - \ldots -mE_8 - (m-1)E_9|.$$
\end{Lem}

\section{Positivité et courbes rationnelles}
\subsection{Un exemple de surface rationnelle dont le diviseur anticanonique n'est pas nef mais s'intersecte positivement avec toute courbe rationnelle}
Soient $9$ points en position très générale dans $\mathbb{P}^2$ de sorte que par ces $9$ points passe une unique cubique lisse $C$. On complète ces neuf points par un dixième toujours sur $C$ et on note $\mu : X \longrightarrow \xe{P}^2$ l'éclatement du plan en ces $10$ points. La transformée stricte $C'$ de $C$ par $\mu$ est un membre irréductible de $|-K_X|$ dont l'intersection avec toute courbe rationnelle est strictement positive. Cependant, comme $-K_X^2 = C'^2=-1$, $-K_X$ n'est pas nef.    

\subsection{Un exemple de surface rationnelle dont le diviseur anticanonique n'est pas pseudoeffectif mais s'intersecte positivement avec toute courbe rationnelle}
Soient $13$ points $x_1,\ldots , x_{13}$ en position très générale dans $\mathbb{P}^2$ de sorte qu'entre autre, passe par ces $13$ points un pinceau de quartiques et aucune cubique. On note $\mu : X \longrightarrow \xe{P}^2$ l'éclatement du plan en ces $13$ points. Puisque la transformée stricte $C$ par $\mu$ de toute quartique passant par ces $13$ points vérifie $-K_X \cdot C = -1$ et que ces courbes couvrent un ouvert dense de $X$, le diviseur anticanonique n'est pas pseudoeffectif. Cependant, le diviseur anticanonique de $X$ s'intersecte positivement avec toute courbe rationnelle, comme le montre le lemme suivant (\cite{GP} lemme 4.2 page 74) : 
\begin{Lem}
\label{GP}
Soit un couple $(d,\alpha)$, où $d$ est un entier strictement positif et $\alpha$ un $r$-uplet $\alpha = (a_1,\ldots ,a_r)$. On note $X_r$ l'éclatement de $\xe{P}^2$ en $r$ points distincts en position très générale, $H$ le tiré en arrière d'un diviseur hyperplan de $\xe{P}^2$, $E_i$ la préimage du point d'éclatement $x_i$. On désigne par $\bar{M}_{(0,0)}(X_r)(d,\alpha)$ l'espace de modules des courbes stables non pointées de genre $0$ et de classe $dH-\sum_{i=1}^r a_i E_i$.

Si $3d - 1 - \sum a_i <0$ alors  $\bar{M}_{(0,0)}(X_r)(d,\alpha)$ est vide.
\end{Lem}
En effet, si $C$ est une courbe rationnelle (irréductible et réduite) s'intersectant négativement avec l'anticanonique de $X$, l'image de $C$ dans $\xe{P}^2$ est alors une courbe rationnelle $\Gamma$ de degré $d$ et de multiplicité $a_i$ aux points $x_i$ vérifiant
$$3d-\sum_{i=1}^{13} a_i \leqslant 0.$$
Or une telle courbe n'existe pas d'après le lemme précédent.

\bigskip
 \noindent A.B. \emph{e-mail: broustet@ujf-grenoble.fr}
 
 \smallskip
 
 \noindent\emph{Institut Fourier, UFR de Math\'ematiques, Universit\'e
   de Grenoble 1, UMR 5582, BP 74, 38402 Saint Martin d'H\`eres, FRANCE}

\begin{thebibliography}{A FINIR}

\bibitem[D1]{Dem1}  J.P. Demailly. \textit{Singular Hermitian metrics on positive line bundles},Complex algebraic varieties, Proc. Conf., Bayreuth/Ger. 1990, Lect. Notes Math. 1507, 87-104. 1992.
\bibitem[D2]{Lnm}  M. Demazure. \textit{Surfaces de Del Pezzo. I-V,} Sémin. sur les singularités des surfaces, Cent. Math. Ec. Polytech., Palaiseau 1976-77, Lect. Notes Math. 777, 21-69. 1980.
\bibitem[F]{Ful} W. Fulton. \textit{Algebraic curves}, Mathematics Lecture Note Series, New York-Amsterdam: W.A. Benjamin, Inc. XIII. 1969.
\bibitem[GLS]{GLS}  G.M. Greuel, C. Lossen et E. Shustin. \textit{Geometry of families of nodal curves on the blown-up projective plane}, Trans. Am. Math. Soc. 350, No.1, 251-274. 1998.
\bibitem[GP]{GP}  L. G\"ottsche et R. Pandharipande. \textit{The quantum cohomology of blow-ups of $\xe{P}^2$ and enumerative geometry}, J. Differ. Geom. 48, No.1, 61-90. 1998.
\bibitem[L]{Laz}  R. Lazarsfeld. \textit{Positivity in algebraic geometry. I. Classical setting: line bundles and linear series}, Ergebnisse der Mathematik und ihrer Grenzgebiete 48, Berlin: Springer. 2004.
\end{thebibliography}
\end{document}